\documentclass[reqno,A4paper]{amsart}
\usepackage{amsmath,amssymb}
\usepackage{color}

\theoremstyle{plain}
\newtheorem{thm}{Theorem}[section]
\newtheorem{lem}[thm]{Lemma}
\newtheorem{cor}[thm]{Corollary}
\newtheorem{prop}[thm]{Proposition}

\newtheorem{question}[thm]{Question}

\newtheorem{defn}[thm]{Definition}

\theoremstyle{definition}
\newtheorem{example}[thm]{Example}
\newtheorem{remark}[thm]{Remark}

\newcommand{\Bd}{\mathrm{Bd}\,}

\newcommand{\Cl}{\mathrm{Cl}\,}

\newcommand{\Sn}{\mathrm{Sn}}
\newcommand{\ind}{\mathrm{ind}\,}
\newcommand{\Ind}{\mathrm{Ind}\,}
\newcommand{\Dm}{\mathrm{Dim}\,}
\newcommand{\ord}{\mathrm{ord}\,}

\begin{document}

\title{On structural numbers of topological spaces.}

\author[Vitalij A.~Chatyrko \and  Alexandre Karassev]{Vitalij A.~Chatyrko* \and  Alexandre Karassev**}

\newcommand{\acr}{\newline\indent}

\address{\llap{*\,}Department of Mathematics\acr
Linkoping University\acr
581 83 Linkoping\acr
SWEDEN}

\email{vitalij.tjatyrko@liu.se}

\address{\llap{**\,}Department of Computer Science and Mathematics\acr
Nipissing University\acr
100 College Drive, Box 5002\acr
North Bay\acr
Ontario P1B 8L7\acr
CANADA}
\email{alexandk@nipissingu.ca}

\subjclass[2010]{Primary 54F45; Secondary 54A10}
\keywords{Structural numbers, small inductive dimension, covering dimension, hereditarily normal spaces.}

\begin{abstract} Zero-dimensional structural numbers $\mathcal Z_0^{\ind}$  and $\mathcal Z_0^{\dim}$ w.r.t. dimensions $\ind$ and  $\mathrm{dim}$   were introduced  by Georgiou, Hattori, Megaritis, and Sereti. Somewhat similarly,  we define structural numbers $\Sn^{\mathcal A}$ for
different subclasses $\mathcal A$ of the class of hereditarily normal $T_1$-spaces. In particular, we show that:
\begin{itemize}
\item[(a)] for any  metrizable space  $X$ with $\dim X = n \geq 0$ we have $1 \leq \Sn^{\mathcal M_{dim}}X \leq n+1$;
\item[(b)] for any  countable-dimensional metrizable space  $Y$  we have $1 \leq \Sn^{\mathcal M_{dim}}Y \leq \aleph_0$,

where $\mathcal M_{dim}$ is the class of metrizable  spaces $Z$ with $\dim Z = 0.$

\end{itemize}

\end{abstract}

\maketitle

\section{Introduction}
For a set $X$, let $\Sigma(X)$ denote the collection of all topologies on $X$, partially ordered by inclusion. Then  $\Sigma(X)$ is a complete lattice in which the meet of a collection of topologies is their intersection, while the join is the topology with their union as a subbasis. This lattice has been the subject of study since it was first defined by Birkhoff in \cite{Bi}, see
\cite{LA} for an overview of some earlier results on this topic.

Let $\tau$ and $\mu$ be topologies on a set $X$. If $\tau \subseteq \mu$ then the topology $\mu$ is called \textit{an extension} of $\tau$. Extension of topologies remains an active area of research (see, for example, \cite{ChH}, \cite{SW}, \cite{MW}).

For a set $X$ and a topology $\tau$ on $X$, let  $\Sigma(X, \tau)$
denote the complete sublattice of $\Sigma(X)$ consisting of all extensions of  $\tau$.
Note that  $\tau$
is the least element of  $\Sigma(X, \tau)$.

Let $\mathcal A$ be a class of topological spaces. If a space $(X, \tau)$ is from the class $\mathcal A$ we will say that the topology $\tau$ is of the class $\mathcal A$. Let $\tau$ and $\mu$ be topologies on a set $X$ and  $\tau \subseteq \mu$. If the topology $\mu$ is of a class $\mathcal A$ then we will say that
$\mu$ is an $\mathcal A$-extension of $\tau$.

In this paper, the case of interest is when the topology $\tau$ is the meet of a subcollection of $\Sigma(X, \tau)$ consisting of $\mathcal A$-extensions of $\tau$. Moreover, we will investigate the interplay between extension of topologies and dimension theory.

The following definition will be used throughout the paper.

\begin{defn} Let $\mathcal A$ be a class of topological spaces and $(X, \tau)$ be a topological space.
\textit{The structural number of $(X, \tau)$ w.r.t. the class $\mathcal A$}, denoted by $\Sn^{\mathcal A}(X, \tau)$, is the smallest cardinal number $\kappa$ for which there exists a family of $\mathcal A$-extensions
$\tau_i, i \in I,$ of $\tau$  such that $|I| = \kappa$ and $\tau = \cap_{i \in I} \tau_i.$
If there is no such family then we define  $\Sn^{\mathcal A}(X, \tau) = \infty$. Here, we  adopt the convention  that $\infty > \kappa$ for any cardinal number $\kappa$.
\end{defn}

It is easy to see that if $\mathcal A, \mathcal B$ are two classes of  spaces such that
$\mathcal A \subseteq \mathcal B$ then $Sn^{\mathcal B}(X, \tau) \leq Sn^{\mathcal A}(X, \tau)$ for any  space $(X, \tau)$. Moreover, if $\mathbb R$ is the set of real numbers, $\tau_d$ is the discrete topology on $\mathbb R$, $\mathcal A_1$ is the class of separable spaces, $\mathcal A_2$ is the class of metrizable spaces then
$Sn^{\mathcal A_1}(X, \tau_d) = \infty$ and $Sn^{\mathcal A_2}(X, \tau_d)  =1$.
However, if $\tau_E$ is the standard Euclidean topology on $\mathbb R$, $\mathcal A_3$ is the class of separable metrizable zero-dimensional spaces then
$Sn^{\mathcal A_3}(X, \tau_E) = 2$ (see Remark \ref{remark_zero}).

\begin{remark} Note that if $\mathcal A$ is the class of zero-dimensional  spaces w.r.t. the small inductive dimension $\ind$
(resp. the covering dimension $\dim$) and  $\Sn^{\mathcal A} (X, \tau) \ne \infty$ then $\Sn^{\mathcal A} (X, \tau)$ is equal to $Z_0^{\ind} (X)$ (resp. $Z_0^{\dim} (X)$), \textit{the zero-dimensional structural number w.r.t. $\ind$ (resp. $\mathrm{dim}$)}, introduced  by
Georgiou, Hattori, Megaritis, and Sereti in \cite{GHMS}. Note also that the main focus of   \cite{GHMS}  is on $T_1$ non-Hausdorff spaces. We, however, consider spaces which are at least normal, and often metrizable.
\end{remark}

This study is devoted to the following question:
\textit{ given subclass of the class of hereditarily normal $T_1$-spaces $\mathcal A$ and a hereditarily normal $T_1$-space $(X, \tau)$
find the structural number $\Sn^{\mathcal A} (X, \tau)$.}

In particular, we show that
\textit{for any  metrizable space  $X$ with $\dim X = n \geq 0$ we have $1 \leq \Sn^{\mathcal M_{dim}}X \leq n+1$,
where $\mathcal M_{dim}$ is the class of metrizable  spaces $Z$ with $\dim Z = 0.$
Moreover, if $\dim (X, \tau) \geq 1$ then
$2 \leq \Sn^{\mathcal M_{\dim}}(X, \tau) \leq n+1$.}

For standard notions we refer to \cite{E1} and \cite{E2}.

\section{A usuful construction}

The following modification of topologies was introduced by Bing \cite{B} and by Hanner \cite{H}.

Let $(X, \tau)$ be a  space and $M \subseteq X$. Define the topology $\tau(M)$
on the set $X$ as follows:
$$\tau(M) = \{U \cup K: U \in \tau\mbox{ and } K \subseteq X\setminus M\}.$$

It is easy to see that   each point of the set $X \setminus M$ is open
in the space $(X, \tau(M))$. The subspace $M$ of $(X, \tau(M))$ is closed and its topology coincides with the topology induced on $M$ by $\tau$. Moreover, $\tau \subseteq  \tau(M)$.

\begin{prop}\label{intersection} Let $(X, \tau)$ be a  space and $M_i \subseteq X, i \in I.$ We have the following:
\begin{itemize}
\item[(1)] if $X = \cup_{i \in I} M_i$ then $\tau = \cap_{i \in I} \tau(M_i)$;
\item[(2)] if $(X, \tau)$ has no isolated points and  $\tau = \cap_{i \in I} \tau(M_i)$ then
$X = \cup_{i \in I} M_i$.
\end{itemize}
\end{prop}

\begin{proof} (1) We need to check only the inclusion $\tau \supseteq \cap_{i \in I} \tau(M_i)$.  Consider any  $\emptyset \ne O \in \cap_{i \in I} \tau(M_i)$ and any $x \in O$. Then there exists   $i_x \in I$ such that $x \in M_{i_x}$. Since $O \in \tau(M_{i_i})$ there exists $U_x \in \tau$ such that $x \in U_x \subseteq O$.
Note that $O = \cup_{x \in O} U_{i_x} \in \tau$.

(2) Assume that the conclusion does not hold, i.e. there exists $x \in X \setminus \cup_{i \in I} M_i$.
By the Bing-Hanner's construction the set $\{x\}$ is open in each topology $\tau(M_i), i \in I$.
Hence,  $\{x\} \in  \cap_{i \in I} \tau(M_i) = \tau.$ We have a contradiction.
\end{proof}

The following simple facts about the  spaces  $(X, \tau(M))$ can be found in \cite[Example 5.1.22]{E1}:
\begin{itemize}
\item[(a)] if $(X, \tau)$ is a $T_i$-space, where $i= 0, 1, 2, 3, 3\frac{1}{2}$, then so is  $(X, \tau(M))$;
\item[(b)] if $(X, \tau)$ is a $T_4$-space  and $M$ is closed then so is $(X, \tau(M))$.
\end{itemize}

Moreover, in the realm of $T_1$-spaces we have the following \cite[Problems 5.5.2]{E1}:
\begin{itemize}
\item[(c)] if $(X, \tau)$ is a hereditarily normal space or a hereditarily collectionwise normal space then
so is $(X, \tau(M))$;
\item[(d)] for a perfectly normal space $(X, \tau)$ the space $(X, \tau(M))$ is perfectly normal if and only if
the set $M$ is a $G_\delta$-set in $(X, \tau)$;
\item[(e)] for a metrizable space $(X, \tau)$ the space $(X, \tau(M))$ is metrizable if and only if
the set $M$ is a $G_\delta$-set in $(X, \tau)$.
\end{itemize}

\section{The Bing-Hanner's constuction and dimensions}

In what follows, all spaces are assumed to be at least $T_1$. As soon as the only topology considered on a set $X$ is $\tau$, we will write $X$ in place of $(X, \tau).$
In this section, we established several facts that connect the Bing-Hanner's construction with various types of dimension. For the reader's convenience,  we also include the basic definitions of dimension theory.

Recall (cf. \cite{AN}) that a subset $L$ of a  space $X$ is a \textit{partition} between two disjoint subsets $A$ and $B$ of $X$ if there exist open disjoint subsets $U, V$ of $X$ such that $A \subseteq U$, $B \subseteq V$ and $L = X \setminus (U \cup V)$.

\textit{The large inductive dimension} $\Ind$ is defined as follows.

Let $X$ be a normal space and let $n$ denote a non-negative integer. Then

\begin{itemize}
\item[(i)] $\Ind X =-1$ if and only if $X = \emptyset$.
\item[(ii)] $\Ind X \leq n$ if for every pair of closed disjoint subsets $A$ and $B$ of $X$ there exists
a partition $C$ between $A$ and $B$ such that $\Ind C \leq n-1$.
\item[(iii)] $\Ind X = n$ if $\Ind X \leq n$ and the inequality $\Ind X \leq n-1$ does not hold.
\item[(iv)] $\Ind X = \infty$ if the inequality $\Ind X \leq n$ does not hold for any $n$.

\end{itemize}

For the proof of the next lemma, see, for example, \cite[Proposition 4.6, the case $n = 0$]{AN}.

\begin{lem}\label{partition} Let $X$ be a hereditarily normal space and $M$ be a subspace of $X$  with $\Ind M =0$. Then for every disjoint closed subsets $A$ and $B$ of $X$ there exists a partition $L$  between $A$ and $B$  in $X$ such that $L \cap M = \emptyset.$
\end{lem}

\begin{prop}\label{zero-subspace} Let $(X, \tau)$ be a hereditarily normal space
(resp. a hereditarily collectionwise normal space) and $M$ be a subspace of $(X, \tau)$  with $\Ind M =0$. Then  $\Ind (X, \tau(M)) = 0$ and the space $(X, \tau(M))$ is hereditarily normal (resp. hereditarily collectionwise normal).

In addition, if the space $(X, \tau)$ is a perfectly normal space (resp. a metrizable space) then
the space  $(X, \tau(M))$ is perfectly normal (resp. metrizable) iff $M$ is a $G_\delta$-set in $(X, \tau)$.
\end{prop}

\begin{proof} Since the space $(X, \tau)$ is hereditarily normal, the space
$(X, \tau(M))$ is also hereditarily normal by the fact $(c)$ from page 3. Hence by Lemma~\ref{partition} for every disjoint closed subsets $A$ and $B$  of $(X, \tau(M))$ there exists a partition $L$ between $A$ and $B$ in $(X, \tau(M))$
 such that $L \cap M = \emptyset$.  Since $L \subseteq X \setminus M$, the subspace $L$ of $(X, \tau(M))$ is open. This easily implies that the empty set is a partition between $A$ and $B$ in the space
$(X, \tau(M))$. We have  that $\Ind (X, \tau(M)) = 0$. The rest follows from the facts $(c), (d)$ and $(e)$, mentioned on page 3.
\end{proof}

In the definition of the covering dimension $\dim$ the concept of the order of a cover is applied.

Recall that \textit{the order} of a family $\nu$ of subsets of a topological space $X$ is defined as follows.

\begin{itemize}
\item[(i)] $\ord \nu = -1$ if $\nu$ consists of the empty set alone.
\item[(ii)] $\ord \nu = n$, where $n$ is a non-negative integer, if the intersection of any $n+2$ distinct elements of $\nu$ is empty and there exists $n+1$ distinct elements of $\nu$ whose intersection is not empty.
\item[(iii)] $\ord \nu = \infty$ if for every positive integer $n$ there exist $n$ distinct elements of $\nu$ whose intersection is not empty.
\end{itemize}

Recall that a non empty family $\mathcal C$ of subsets of $X$ is called \textit{a cover} of $X$ if the union of all elements of $\mathcal C$ is equal to $X$. Moreover, a cover $\mathcal C$ of $X$ is called \textit{open}
if all elements of $\mathcal C$ are open subsets of $X$. Furthermore, a family $\nu$ of subsets of $X$ is called \textit{a refinement} of a family $\mathcal C$ of subsets of $X$, written $\nu > \mathcal C$, if every element of $\nu$ is contained in an element of $\mathcal C$.

\textit{The covering dimension} $\dim$ is defined as follows.

Let $X$ be a normal space and let $n$ denote an integer $\geq -1$. Then

\begin{itemize}
\item[(i)] $\dim X =-1$ if and only if $X = \emptyset$.
\item[(ii)] $\dim X \leq n$ if for every finite open cover $\mathcal C$ of $X$  there exists
a finite open cover $\nu$ of $X$ such that $\nu > \mathcal C$ and  $\ord \nu \leq n$.
\item[(iii)] $\dim X = n$ if $\dim X \leq n$ and the inequality $\dim X \leq n-1$ does not hold.
\item[(iv)] $\dim X = \infty$ if the inequality $\dim X \leq n$ does not hold for any $n$.

\end{itemize}

Recall \cite[Theorem 3.1.28]{E2} that for every normal $T_1$-space $X$ we have $\dim X \leq \Ind X$, and   the conditions $\Ind X =0$ and $\dim X = 0$ are equivalent (see \cite[Theorem 1.6.10]{E2}).

\begin{example}
 Let $(X, \tau)$ be a separable metrizable one-dimensional space.
Let also $\mathcal B = \{B_i\}_{i=1}^\infty$ be a countable base for $(X, \tau)$
such that $\dim Bd B_i \leq 0$ for each $i$. Put $A = \cup_{i=1}^\infty Bd B_i$ and $B = X \setminus A$.
Note that $\dim A = \dim B =0$, and the set $B$ is a non-countable $G_\delta$-set but the set $A$ is not a $G_\delta$-set. Then $\tau = \tau(A) \cap \tau(B)$ by Proposition \ref{intersection}. Furthermore, by Proposition \ref{zero-subspace}
we have  $\dim (X, \tau(A)) = \dim (X, \tau(B)) = 0$,  the space $(X, \tau(B))$ is  metrizable
(and non-separable if $A$ is not countable) and the space $(X, \tau(A))$
is hereditarily collectionwise normal, non-metrizable and non-separable.
\end{example}

\begin{example} Consider the real line $(\mathbb R, \tau_E)$, where $\mathbb R$ is the set of real numbers and $\tau_E$ is the Euclidean topology. Note that $\mathbb R = \mathbb Q \cup \mathbb P$, where $\mathbb Q$ (resp. $\mathbb P$) is the set of rational (resp. irrational) numbers. Then $\tau_E = \tau_E(\mathbb Q) \cap \tau_E(\mathbb P)$ by Proposition \ref{intersection}. Furthermore, by Proposition \ref{zero-subspace}
we have that $\dim (\mathbb R, \tau_E(\mathbb P)) = \dim (\mathbb R, \tau_E(\mathbb Q)) = 0$,  the space $(\mathbb R, \tau_E(\mathbb P))$ is separable metrizable and the space $(\mathbb R, \tau_E(\mathbb Q))$
is hereditarily collectionwise normal,  non-metrizable and non-separable.
\end{example}

\textit{The small inductive dimension} $\ind$ is defined as follows.

Let $X$ be a regular space and let $n$ denote a non-negative integer. Then

\begin{itemize}
\item[(i)] $\ind X =-1$ if and only if $X = \emptyset$.
\item[(ii)] $\ind X \leq n$ if for  any closed subset $A$ of $X$ any point $x \in X \setminus A$ and there exists
a partition $C$ such that $\ind C \leq n-1$.
\item[(iii)] $\ind X = n$ if $\ind X \leq n$ and the inequality $\ind X \leq n-1$ does not hold.
\item[(iv)] $\ind X = \infty$ if the inequality $\ind X \leq n$ does not hold for any $n$.

\end{itemize}

Recall (\cite[Theorem 1.6.3]{E2}) that for every normal space $X$ we have $\ind X \leq \Ind X$.

\begin{remark}\label{Roy}
Note  that there exists a completely metrizable space $R$ such that $\ind R = 0$ and $\Ind R = 1$  \cite{R}.
\end{remark}

\begin{lem}\label{partition_ind} Let $X$ be a hereditarily normal space and $M$ be a subspace of $X$  with $\ind M =0$. Then for every point $x \in M$ and any closed subset $A$ of $X$ such that $x\notin A$  there exists a partition $L$ between $x$ and $A$ in $X$  such that $L \cap M = \emptyset.$
\end{lem}

\begin{proof} Our proof is somewhat similar to the proof of Lemma \ref{partition}. Consider an open subset $O$ of $X$ such that $A \subseteq O$ and $x \notin \Cl_X O$. If $\Cl_X O \cap M = \emptyset$, consider as $L$ any partition between $A$ and $X \setminus O$ in $X$. Otherwise, the empty set is a partition between the point $x$ and $\Cl_X O \cap M$ in the subspace $M$ of $X$, i.e. there are disjoint closed in $M$ subsets $E$ and $F$ such that $x \in E$, $\Cl_X O \cap M \subseteq F,$ and $E \cup F = M.$ Note that the sets $E$ and
$F \cup A$ are separated in $X$, i.e. $\Cl_XE \cap (F \cup A) = E \cap \Cl_X (F \cup A) = \emptyset$. Since $X$ is hereditarily normal, there exist disjoint open subsets $U$ and $V$ of $X$ such that $E \subseteq U$ and $F \cup A \subseteq V$. Put $L = X \setminus (U \cup V)$ and note that $L \cap M = \emptyset$.
\end{proof}

\begin{remark} Note that  the condition $x \in M$ is essential in Lemma~\ref{partition_ind} (see \cite[Problem 4.1.C]{E2}).
\end{remark}

Similarly to Proposition \ref{zero-subspace} we can prove the following.

\begin{prop}\label{zero-subspace_ind} Let $(X, \tau)$ be a hereditarily normal space
(resp. a hereditarily collectionwise normal space) and $M$ be a subspace of $(X, \tau)$  with $\ind M =0$. Then  $\ind (X, \tau(M)) = 0$ and the space $(X, \tau(M))$ is hereditarily normal (resp. hereditarily collectionwise normal).

In addition, if the space $(X, \tau)$ is a perfectly normal space (resp. a metrizable space) then
the space  $(X, \tau(M))$ is perfectly normal (resp. metrizable) if and only if $M$ is a $G_\delta$-set in $(X, \tau)$.
\end{prop}

\begin{proof} Since the space $(X, \tau)$ is hereditarily normal, the space
$(X, \tau(M))$ is also hereditarily normal by the fact $(c)$ from page 3. Note that all points from the set $X \setminus M$ are clopen subsets of $(X, \tau(M))$ so we need to check the small inductive dimension at the points from $M$. Consider any point $x \in M$ and any open neighbourhood $O$ of $x$ in $(X, \tau(M))$.  Since the set  $X \setminus O$ is closed in $(X, \tau(M))$  and $x \notin X \setminus O,$ by Lemma \ref{partition_ind} there exists
a partition $L$ between $x$ and $X \setminus O$ in $(X, \tau(M))$ such that $L \cap M = \emptyset$, i.e.  there exist open disjoint subsets $U$ and $V$ of $(X, \tau(M))$ such that $x \in U$, $X \setminus O \subseteq V,$ and $L = X \setminus (U \cup V)$. Note that $\Bd_{(X, \tau(M))} U  \subseteq L$. Consequently,  $\Bd_{(X, \tau(M))} U  \cap M = \emptyset$. This implies that $\Bd_{(X, \tau(M))} U$ is open in  $(X, \tau(M))$, i.e.
$\Bd_{(X, \tau(M))} U  = \emptyset.$ Hence  $\ind_x (X, \tau(M)) = 0.$
\end{proof}

\section{Main results}
Let $\mathcal A$ be a subclass of the class of hereditarily normal spaces, and $\Dm$ be a dimension function.
Define $$\mathcal A_\Dm = \{X \in \mathcal A\colon \Dm(X) = 0\}.$$

\begin{remark} \label{ind_Ind} Recall that for any normal $T_1$-space $X$ we have $\ind X \leq \Ind X$. Hence
$\mathcal A_{Ind} \subseteq \mathcal A_{ind}$ for any subclass $\mathcal A$ of the class of hereditarily normal spaces. In particular, $\Sn^{\mathcal A_{ind}} X \leq \Sn^{\mathcal A_{dim}}X$ for any hereditarily normal $T_1$-space $X$.
\end{remark}

\subsection{Finite-dimensional metrizable case}
\hfill\\
Recall (\cite[Theorem 4.1.3]{E2}) that for every metrizable space $X$ we have $\Ind X = \dim X$.
Therefore in this section we consider the dimension $\dim$ instead of $\Ind$.

The following lemma is a consequence of Theorem 4.1.19 from \cite{E2}.

\begin{lem}\label{G_delta} Let $X$ be a metrizable space and $M$ be a subspace of  $X$ with $\dim M = 0$. Then there exists a $G_\delta$-set $M^*$
in  $X$ such that $M \subseteq M^*$ and $\dim M^* = 0$.
\end{lem}

The following fact is also standard (see \cite[Theorem 4.1.17]{E2}).
\begin{lem}\label{decomposition} A  metrizable space $X$ satisfies the inequality $\dim X \leq n \geq 0$ iff $X$ can be represented as a union of  $n+1$ subspaces $M_1, \dots, M_{n+1}$ such that $\dim M_i \leq 0$ for $i=1, \dots, n+1.$
\end{lem}

Let $\mathcal M$ be the class of metrizable  spaces $X$.

\begin{thm} \label{finite-theorem} Let $(X, \tau)$ be a  metrizable space  with $\dim (X, \tau) = n \geq 0$.
Then there exist subsets $M_1, \dots, M_{n+1}$ of $(X, \tau)$ such that
\begin{itemize}
\item[(i)] $\tau = \cap_{i=1}^{n+1} \tau(M_i),$ and
\item[(ii)]  each space
$(X, \tau(M_i)), i=1, \dots, n+1,$ is metrizable, and

$\dim (X, \tau(M_i)) = 0.$
\end{itemize}

In particular,  $1 \leq \Sn^{\mathcal M_{\dim}}(X, \tau) \leq n+1$.
Moreover, if $\dim (X, \tau) \geq 1$ then
$2 \leq \Sn^{\mathcal M_{\dim}}(X, \tau) \leq n+1$.
\end{thm}

\begin{proof} (i) By Lemmas ~\ref{G_delta} and~\ref{decomposition} there exist
 $G_\delta$-sets  $M_1, \dots, M_{n+1}$ in $(X,\tau)$ such that $\dim M_i \leq 0$ for $i=1, \dots, n+1,$ and $X = \cup_{i=1}^{n+1} M_i$. So
by Proposition \ref{intersection} we have $\tau = \cap_{i=1}^{n+1} \tau(M_i).$

(ii) By Proposition \ref{zero-subspace} each space
$(X, \tau(M_i)), i=1, \dots, n+1,$ is metrizable, and $\dim (X, \tau(M_i)) = 0.$

Hence, $1 \leq \Sn^{\mathcal M_{\dim}}(X, \tau) \leq n+1$.
It is evident that if $\dim (X, \tau) \geq 1$ then
$2 \leq \Sn^{\mathcal M_{\dim}}(X, \tau) \leq n+1$.
\end{proof}

\begin{cor} Let $(X, \tau)$ be a separable metrizable one-dimensional space (for example, the real line
$(\mathbb R, \tau_E)$, where $\mathbb R$ is the set of real numbers and $\tau_E$ is the Euclidean topology). Then
$\Sn^{\mathcal M_{ind}}(X, \tau) = \Sn^{\mathcal M_{dim}}(X, \tau) = 2$.
\end{cor}

\begin{remark} In \cite{GHMS} the authors proved that
$Z_0^{\ind}(\mathbb R, \tau_E) = 2$. Note that
$Z_0^{\ind}(X, \tau) \leq \Sn^{\mathcal M_{ind}}(X, \tau)$
for any  space $(X, \tau)$.

\end{remark}

\begin{cor} $\Sn^{\mathcal M_{ind}}R = 1 < \Sn^{\mathcal M_{dim}}R = 2$, where $R$ is the space from
Remark \ref{Roy}. $\Box$
\end{cor}

\begin{example} Recall that applying the existence of a metrizable space with the properties from Remark \ref{Roy}
van Douwen \cite{D} and Przymusinski \cite{P} constructed a metrizable space $Y$ with $\ind Y = \dim Y = 1$
which can be represented as the union of two closed subspaces $Y_1$ and $Y_2$ such that $\ind Y_1 = \ind Y_2 = 0$. Moreover, $Y$ contains a point $a$ such that $\ind (Y \setminus \{a\}) = 0.$
It is evident that $\Sn^{\mathcal M_{ind}}Y = \Sn^{\mathcal M_{dim}}Y = 2.$
\end{example}

\begin{question} Let $(X, \tau)$ be a  metrizable space  with $\dim (X, \tau) = n \geq 1$.
Does the equality $\Sn^{\mathcal M_{dim}}(X, \tau) = n+1$ hold? In particular, is the equality
$\Sn^{\mathcal M_{dim}}(\mathbb R^2, \tau_E) = 3$ (even, $\Sn^{\mathcal M_{ind}}(\mathbb R^2, \tau_E) = 3$) valid, where $\tau_E$ is the Euclidean topology on
$\mathbb R^2$?
\end{question}

\begin{thm}\label{product} Let $(X, \tau)$ and $(Y, \sigma)$ be metrizable spaces with $\ind (X, \tau) = 0$
and $\dim (Y, \sigma) = n, n \geq 1,$ and $\mu$ be the product of topologies $\tau$ and $\sigma$. Then

$\ind (X \times Y, \mu) \leq n$ and  $2 \leq \Sn^{\mathcal M_{ind}} (X \times Y, \mu) \leq n+1$.
\end{thm}
\begin{proof} The equality $\ind (X \times Y, \mu) \leq n$ is evident.
Let us show that
$\Sn^{\mathcal M_{ind}} (X \times Y, \mu) \leq n+1$. For that as in the proof of Theorem \ref{finite-theorem} we can consider  $G_\delta$-sets  $M_1, \dots, M_{n+1}$ in $(Y,\sigma)$ such that $\dim M_i = 0$ for $i=1, \dots, n+1,$ and $Y= \cup_{i=1}^{n+1} M_i$. Note that for each
$i=1, \dots, n+1,$
the set $N_i = X \times M_i$ is $G_\delta$ in $(X \times Y, \mu)$ and $\ind N_i = 0$. Moreover,
 $X \times Y = \cup_{i=1}^{n+1} N_i$. So by Proposition \ref{intersection} we have $\mu= \cap_{i=1}^{n+1} \mu(N_i).$ Then by Proposition \ref{zero-subspace_ind} for each $i=1, \dots, n+1,$ the space
$(X \times Y, \mu(N_i))$ is metrizable, and $\ind (X \times Y, \mu(N_i)) = 0.$ Hence, $\Sn^{\mathcal M_{\ind}}(X \times Y, \mu) \leq n+1.$
\end{proof}

\begin{cor} For every $n \geq 1$ we have  $\ind (R \times I^n) = n$ and $2 \leq \Sn^{\mathcal M_{ind}} (R \times I^n) \leq n+1$, where $R$ is the space from
Remark \ref{Roy} and $I = [0,1]$ is the closed interval with the Euclidean topology.

In particular, $\Sn^{\mathcal M_{ind}} (R \times I) = 2$.
\end{cor}

Note that $\dim (R \times I^n) = n+1, n \geq 1$ \cite[Theorem 2]{Morita}.

\begin{question}
Do the equalities

$\Sn^{\mathcal M_{dim}} (R \times I^n) = n+2$ and $\Sn^{\mathcal M_{ind}} (R \times I^n) = n+1$ hold for each $n \geq 1$? In particular, is the equality
$\Sn^{\mathcal M_{dim}}(R \times I) = 3$ valid?
\end{question}

\begin{lem}\label{without_iso} Let $(X, \tau)$ be a hereditarily normal space without isolated points, and
$\dim (X, \tau) = n,$ $n \geq 1$. Then the topology $\tau$ cannot be represented as the intersection of topologies $\tau(M_i),$ $i =1, \dots, k,$ for any $k \leq  n$ and any zero-dimensional in the sense of
$\dim$ subsets $M_i,$ $i = 1, \dots, k,$ of $(X, \tau)$.
\end{lem}

\begin{proof} Indeed, if $\tau = \cap_{i=1}^k \tau(M_i)$ for some $k \leq n$ and some
zero-dimensional in the sense of $\dim$ subsets $M_i,$  $i = 1, \dots, k,$ of $(X, \tau)$ then
by Proposition \ref{intersection} we have $X = \cup_{i=1}^k M_i.$ Furthermore, by the addition theorem for $\dim$ (see \cite[3.1.17]{E2}) we get $\dim (X, \tau) \leq k-1 < n$. This is a contradiction.
\end{proof}

\begin{thm}\label{thm_equivalence} A metrizable space $(X, \tau)$ without isolated points satisfies the
inequality $\dim (X, \tau) \leq n,$ $n \geq 1$, if and only if the space  $(X, \tau)$ contains $(n+1)$
subsets $M_1, \dots, M_{n+1}$
such that
\begin{itemize}
\item[(i)] $\tau = \cap_{i=1}^{n+1} \tau(M_i),$ and
\item[(ii)]  each space
$(X, \tau(M_i)), i=1, \dots, n+1,$ is metrizable, and

$\dim (X, \tau(M_i)) = 0.$
\end{itemize}

\end{thm}

\begin{proof}The necessity follows from Theorem \ref{finite-theorem}. The sufficiency follows  Lemma~\ref{without_iso} and the observation that the subsets $M_1, \dots, M_{n+1}$ of $(X, \tau)$ are zero-dimensional in the sense of $\dim$.
\end{proof}

\begin{cor} If a metrizable space $(X, \tau)$ without isolated points satisfies the
equality $\dim (X, \tau) = \infty$ then  for each $n \geq 1$
the space  $(X, \tau)$ contains no
subsets $M_1, \dots, M_{n}$
such that
\begin{itemize}
\item[(i)] $\tau = \cap_{i=1}^{n} \tau(M_i),$ and
\item[(ii)]  each space
$(X, \tau(M_i)), i=1, \dots, n,$ is metrizable, and

$\dim (X, \tau(M_i)) = 0.$
\end{itemize}
\end{cor}

Let $\mathcal SM$ be the class of separable metrizable  spaces. Note that
$\mathcal SM_{ind} = \mathcal SM_{dim}$. Hence, $\Sn^{\mathcal SM_{ind}} X = \Sn^{\mathcal SM_{dim}} X$
for any normal $T_1$-space $X$.

The following is evident.
\begin{prop} Let $\tau$ and $\mu$ be topologies on a set $X$ such that $\tau \subseteq \mu$.
Then $d(X, \tau) \leq d(X, \mu)$, where $d(Z)$ is the density of a space $Z$.
In particular, if the space $(X, \mu)$ is separable then the space $(X, \tau)$ is also separable. $\Box$
\end{prop}

\begin{cor} For any non-separable metrizable space $X$ we have
$\Sn^{\mathcal SM_{dim}} X  = \infty$.
\end{cor}

\begin{remark}\label{remark_zero}  Note that  $\Sn^{\mathcal SM_{dim}}(\mathbb R, \tau_E)  = 2$. In fact,
we can represent $\mathbb R$ as a union of two $G_\delta$-sets $M_1$ and $M_2$ in  $(\mathbb R, \tau_E)$ with $\dim M_i =0.$ For that consider two disjoint countable dense subsets $Q_1$ and $Q_2$ of $(\mathbb R, \tau_E)$, and then put $M_1 = \mathbb R \setminus Q_1$ and $M_2 = \mathbb R \setminus Q_2$.   Observe that $\tau_E = \tau(M_1) \cap \tau(M_2)$ and the spaces $(\mathbb R, \tau(M_i)),$  $i =1,2,$ are separable metrizable and zero-dimensional.
\end{remark}

\begin{prop} Let $(X, \tau)$ be a  separable metrizable space  with $\dim (X, \tau) = n,$ $n \geq 2$.
Then for any zero-dimensional $G_\delta$-subset $M$ of $(X, \tau)$ the metrizable space $(X, \tau(M))$ is not separable.
\end{prop}

\begin{proof} Indeed if $|X \setminus M| \leq \aleph_0$ then $\dim (X \setminus M, \tau|_{(X \setminus M)}) = 0$. Hence, $\dim (X, \tau) \leq 1$. We have a contradiction. So $|X \setminus M| > \aleph_0$ and the
metrizable space $(X, \tau(M))$ is not separable.
\end{proof}

\begin{cor}  For a separable metrizable space $X$ with $\dim(X, \tau) \geq 2$ we cannot replace  ``metrizable" with 	``separable metrizable" in Theorem \ref{thm_equivalence}.(ii).
\end{cor}

\begin{remark} Note that for each integer $n \geq 1$ there is a continuous bijection of the space $\mathbb P$ of irrational numbers onto the Euclidean space ($\mathbb R^n, \tau_E)$, i.e. there is an extension $\nu$ of $\tau_E$ such that the space $(\mathbb R^n, \nu)$ is homeomorphic to $\mathbb P$.
\end{remark}

\begin{question} Let $(X, \tau)$ be a  separable metrizable space  with $\dim (X, \tau) = n,$ $n \geq 2$.
Does the inequality $\Sn^{\mathcal SM_{ind}}(X, \tau) \leq n+1$ hold? In particular, is the inequality
$\Sn^{\mathcal SM_{ind}}(\mathbb R^2, \tau_E) \leq 3$ valid, where $\tau_E$ is the Euclidean topology on
$\mathbb R^2$?
\end{question}

\subsection{Infinite-dimensional metrizable case}
\hfill\\
Recall \cite[Proposition 5.1.3]{E2} that a metrizable space is \textit{countable-dimensional (in the sense of $\dim$)} if
$X$ can be represented  as a union of a sequence  of subspaces $X_1, X_2, \dots$ such that $\dim X_i \leq 0$
for $i =1, 2, \dots$.

Similarly to Theorem \ref{finite-theorem}, using Propositions~\ref{intersection} and~\ref{zero-subspace} one can prove the following.
\begin{thm} Let $(X, \tau)$ be a  metrizable countable-dimensional  space.
 Then $1 \leq \Sn^{\mathcal M_{dim}}(X, \tau) \leq \aleph_0$.
\end{thm}

\begin{question} Let $(X, \tau)$ be a  metrizable countable-dimensional  space with $\dim X = \infty$.
Does the equality $\Sn^{\mathcal M_{dim}}(X, \tau) = \aleph_0$ hold?
\end{question}

In general case we have the following.
\begin{prop} Let $(X, \tau)$ be a  metrizable  space.
Then

$\Sn^{\mathcal M_{dim}}(X, \tau) \leq |X|$. $\Box$
\end{prop}

\begin{proof} Note that $X = \cup_{x \in X} \{x\}$.
Then by Propositions~\ref{intersection} and  ~\ref{zero-subspace}  we get the inequality.
\end{proof}

\subsection{Non-metrizable case}
Let $\mathcal HN$  (resp. $\mathcal HCN$ ) be the class of hereditarily normal (resp. collectionwise normal)  spaces $X$.

Using the facts $(c)$ and  $(d),$ mentioned on page 3, and Propositions~\ref{intersection},~\ref{zero-subspace}, or~\ref{zero-subspace_ind}, we get

\begin{thm} Let $(X, \tau)$ be a hereditarily normal (resp. collectionwise normal) space
such that $X = \cup_{i \in I} X_i$, where each $X_i$ is a subspace of $(X, \tau)$ with $\Dm X_i =0,$ where $\Dm$ is either $\dim$ or $\ind$.
Then $1 \leq \Sn^{\mathcal HN_{\Dm}}(X, \tau) \leq |I|$ (resp. $1 \leq \Sn^{\mathcal HCN_{\Dm}}(X, \tau) \leq |I|$).
\end{thm}

\end{document}